\documentclass[12pt,a4paper]{article}
\usepackage[cp1251]{inputenc} 
\usepackage{amsfonts}

\newcommand{\di}{\displaystyle}

\newcommand{\al}{\alpha}
\newcommand{\be}{\beta}
\newcommand{\ga}{\gamma}

\newcommand{\ee}{\varepsilon}

\newcommand{\iy}{\infty}

\begin{document}

\begin{center}
{\large\bf
ASYMPTOTICS OF SOLUTIONS OF DIFFERENTIAL EQUATIONS WITH A SPECTRAL
PARAMETER}\\[0.3cm]
{\bf V.\,A. Yurko} \\[0.5cm]
\end{center}

\vspace{0.5cm}

{\bf Abstract.} The main goal of this paper is to construct the so-called Birkhoff-type
solutions for linear ordinary differential equations with a spectral
parameter. Such solutions play an important role in direct and inverse problems of spectral theory.
In Section 1, we construct the Birkhoff-type
solutions for $n$-th order differential equations. Section 2 is devoted to
first-order systems of differential equations.

\medskip

{\bf Keywords:} Birkhoff-type solutions; linear ODU; higher-order differential operators; first-order differential systems;
asymptotics of solutions.

\medskip

{\bf AMS Mathematics Subject Classification (2020):} 34A30 34D05 34E10

\vspace{1cm}

The main goal of this paper is to construct the so-called Birkhoff-type
solutions for linear ordinary differential equations with a spectral
parameter. Such solutions play an important role in many problems of the
spectral theory (see, for example, [1] and references therein). Moreover,
they also appear in the inverse problem theory ([2]-[3]).

This paper contains two sections. In Section 1 we construct the Birkhoff-type
solutions for n-th order differential equations, and Section 2 is devoted to
first-order systems of differential equations. Other results related to this area one can find in [4]-[10].

\bigskip

\begin{center}
{\bf I. Birkhoff-type solutions for arbitrary order differential
equations}\\[0.3cm]
\end{center}

{\bf 1.1.}$\;$ In this section we study the differential equation
of order $n\ge 2$:
$$
y^{(n)}+\di\sum_{m=0}^{n-2}p_m(x)y^{(m)}=\rho^{n}y,
\quad 0\le x\le T\le \iy                                      \eqno(1.1)
$$
on the finite interval $(T< \iy)$ or on the half-line $(T=\iy)$.
Here $\rho$ is the spectral parameter, and $p_m(x)\in L(0,T)$
are complex-valued integrable functions.

\bigskip
Our goal here is to construct a fundamental system of solutions (FSS)
$\{y_k(x,\rho)\}_{k=\overline{1,n}}$ of equation (1.1) such that
$$
y_k(x,\rho) \sim \exp(\rho R_k x), \quad |\rho|\to\iy,
$$
where $R_1,R_2,\ldots,R_n$ are the roots of the equation $R^{n}-1=0.$
We note that the functions $\{\exp(\rho R_k x)\}_{k=\overline{1,n}}$
form the FSS for the "simplest" equation $y^{(n)}=\rho^{n}y,$ when
in (1.1) $p_m(x)=0,\; m=\overline{0,n-2}.$
\bigskip

It is easy to see that the $\rho$-plane can be partitioned into sectors $S$
of angle $\di\frac{\pi}{n}\quad \Big(\arg\rho \in \Big(\di\frac{\mu\pi}{n},
\di\frac{(\mu+1)\pi}{n}\Big),\; \mu=\overline{0,2n-1}\Big)$ in which the
roots $R_1, R_2,\ldots , R_n$ of the equation $R^{n}-1=0$ can be numbered
in such a way that
$$
Re(\rho R_1)< Re(\rho R_2)<\ldots<Re(\rho R_n),\quad \rho\in S. \eqno(1.2)
$$

Fix $\al \in [0,T),\; k=\overline{1,n},$ and a sector $S$ with the
property (1.2). Consider the following integro-differential equation
with respect to $y_k(x,\rho),\; x \in [\al,T],\; \rho\in\overline{S}$:
$$
y_k(x,\rho)=\exp(\rho R_k x)- \di\frac{1}{n\rho^{n-1}}
\di\int_\al^{x}\Big(\di\sum_{j=1}^{k}R_j\exp(\rho R_j (x-t))\Big)
\Big(\di\sum_{m=0}^{n-2} p_m(t) y_k^{(m)}(t,\rho)\Big)\,dt
$$
$$
+\di\frac{1}{n\rho^{n-1}}
\di\int_x^{T}\Big(\di\sum_{j=k+1}^{n} R_j\exp(\rho R_j(x-t))\Big)
\Big(\di\sum_{m=0}^{n-2}p_m(t) y_k^{(m)}(t,\rho)\Big)\,dt.    \eqno(1.3)
$$

\bigskip
{\bf Remark 1.1.} We note that the general solution of the
differential equation
$$
y^{(n)}=\rho^{n}y+f(x), \quad 0\le x\le T,
$$
has the form
$$
y(x,\rho)=
\di\sum_{j=1}^{n}C_j\exp(\rho R_j x)+
\di\frac{1}{n\rho^{n-1}}\di\sum_{j=1}^{n}
\di\int_{\ga_j}^{x} R_j\exp(\rho R_j(x-t)) f(t)\,dt,
$$
where $\gamma_j\in [0,T]$ are arbitrary fixed numbers.
Clearly, (1.3) corresponds to the case
$$
C_j=
\left\{ \begin{array}{ll}
0,\quad & j\ne k, \\
1,\quad & j=k,
\end{array}\right. \qquad
\ga_j=
\left\{ \begin{array}{ll}
\al,\quad & j\le k,\\
T,\quad & j>k.
\end{array}\right.
$$

\bigskip
Let us now transform (1.3) to a system of linear integral equations.
Differentiating (1.3) with respect to $x$, we get for
$\nu=\overline{0,n-1},$
$$
y_k^{(\nu)}(x,\rho)= (\rho R_k)^{\nu}\exp(\rho R_k x)
$$
$$
-\di\frac{1}{n\rho^{n-1}}\di\int_\al^{x}
\Big(\di\sum_{j=1}^{k} R_j (\rho R_j)^{\nu}\exp(\rho R_j(x-t))\Big)
\Big(\di\sum_{m=0}^{n-2} p_m(t) y_k^{(m)}(t,\rho)\Big)\,dt+
$$
$$
\di\frac{1}{n\rho^{n-1}}\di\int_x^{T}
\Big(\di\sum_{j=k+1}^{n}R_j(\rho R_j)^{\nu}\exp(\rho R_j(x-t))\Big)
\Big(\di\sum_{m=0}^{n-2} p_m(t) y_k^{(m)}(t,\rho)\Big)\,dt.    \eqno(1.4)
$$
Denote
$$
z_{\nu k}(x,\rho)=(\rho R_k)^{-\nu}\exp(-\rho R_k x)y_k^{(\nu)}(x,\rho).
$$
Then
$$
y_k^{(\nu)}(x,\rho)= (\rho R_k)^{\nu}\exp(\rho R_k x)z_{\nu k}(x,\rho),
$$
and (1.4) implies
$$
z_{\nu k}(x,\rho)=
$$
$$
1-\di\frac{1}{n\rho^{n-1}}\di\int_\al^{x}
\Big(\di\sum_{j=1}^{k} R_j^{\nu +1} R_k^{-\nu}\exp(\rho (R_j-R_k)(x-t))\Big)
\Big(\di\sum_{m=0}^{n-2} p_m(t) (\rho R_k)^{m}z_{m k}(t,\rho)\Big)\,dt
$$
$$
+\di\frac{1}{n\rho^{n-1}}\di\int_x^{T}
\Big(\di\sum_{j=k+1}^{n} R_j^{\nu +1} R_k^{-\nu}
\exp(\rho (R_j-R_k)(x-t))\Big)\Big(\di\sum_{m=0}^{n-2}
p_m(t)(\rho R_k)^{m}z_{m k}(t,\rho)\Big)\,dt,                 \eqno(1.5)
$$
or
$$
z_{\nu k}(x,\rho)=1+\di\sum_{m=0}^{n-2}
\int_\al^{T}A_{\nu m k}(x,t,\rho)z_{m k}(t,\rho)\,dt,
\quad \nu=\overline{0,n-1},                                   \eqno(1.6)
$$
where
$$
A_{\nu m k}(x,t,\rho)=
\left\{\begin{array}{ll}
-\di\frac{p_m(t)}{n\rho^{n-1-m}}
\di\sum_{j=1}^{k} R_j^{\nu +1} R_k^{m-\nu}\exp(\rho (R_j-R_k)(x-t)),
\quad & x\ge t, \\
\di\frac{p_m(t)}{n\rho^{n-1-m}}
\di\sum_{j=k+1}^{n} R_j^{\nu +1} R_k^{m-\nu}\exp(\rho (R_j-R_k)(x-t)),
\quad & x< t.
\end{array}\right.                                            \eqno(1.7)
$$

For fixed $\rho\in \overline{S}$ and $k=\overline{1,n},$ we consider
(1.6) as a system  of linear integral equations with respect to
$z_{\nu k}(x,\rho),\; x\in [\al,T].$ By virtue of (1.7) and (1.2) we have
$$
\di\int_\al^{T}|A_{\nu m k}(x,t,\rho)|\,dt\;\le
\di\frac{1}{n|\rho|^{n-1-m}}\di\int_\al^{T}|p_m(t)|\,dt.      \eqno(1.8)
$$
Denote
$$
\rho_\al:=\max_{m=\overline{0,n-2}}
\Big(2\di\int_\al^{T}|p_m(t)|\,dt\Big)^{\frac{1}{n-1-m}}.
$$
It follows from (1.8) that for $|\rho|\ge \rho_\al,\; \rho \in\overline{S},
\;x \in [\al,T],\; k=\overline{1,n},\; \nu=\overline{0,n-1}$:
$$
\di\sum_{m=0}^{n-2}\di\max_{\al\le x\le T}
\di\int_\al^{T}|A_{\nu m k}(x,t,\rho)|\,dt\;\le \di\frac{1}{2}. \eqno(1.9)
$$
Solving (1.6) by the method of successive approximations and using (1.9), we
obtain that for $|\rho|\ge \rho_\al,\; \rho\in \overline{S},\; x\in [\al,T],
\;k=\overline{1,n},$ system (1.6) has a unique solution $z_{\nu k}(x,\rho),
\;\nu=\overline{0,n-1},$ and $|z_{\nu k}(x,\rho)|\le 2$. Substituting this
estimate into the right-hand side of (1.6) and using (1.8), we get
$$
|z_{\nu k}(x,\rho)-1|\le\di\frac{C}{|\rho|},\quad |\rho|\ge \rho_\al,\;
\rho\in\overline{S},\; x\in [\al,T],\; k=\overline{1,n},\;
\nu=\overline{0,n-1}.
$$
In other words, for $|\rho|\to\iy,\; \rho \in\overline{S},$
$$
z_{\nu k}(x,\rho)=1+O\Big(\di\frac{1}{\rho}\Big),\quad
k=\overline{1,n},\; \nu=\overline{0,n-1},                   \eqno(1.10)
$$
uniformly in $x \in [\al,T].$ Moreover, one can obain more precise
asymptotic formulae for the functions $z_{\nu k}(x,\rho)$ than
(1.10). Indeed, substituting (1.10) into the right-hand side of
(1.5), we calculate
$$
z_{\nu k}(x,\rho)=1-\di\frac{1}{n\rho}\di\int_\al^{x}
\Big(\di\sum_{j=1}^{k} R_j^{\nu +1} R_k^{n-\nu-2}
\exp(\rho (R_j-R_k)(x-t))\Big)p_{n-2}(t)\,dt
$$
$$
+\di\frac{1}{n\rho}\di\int_x^{T}\Big(\di\sum_{j=k+1}^{n}
R_j^{\nu +1} R_k^{n-\nu-2}\exp(\rho (R_j-R_k)(x-t))\Big)
p_{n-2}(t)\,dt+O\Big(\di\frac{1}{\rho^{2}}\Big).
$$
The terms with $j\ne k$ give us $o\Big(\di\frac{1}{\rho}\Big)$ as
$|\rho|\to\iy,\;\rho\in \overline{S},$ uniformly in $x \in [\al,T].$ Hence
$$
z_{\nu k}(x,\rho)=1+\di\frac{\be_1(x)}{\rho R_k}+
o\Big(\di\frac{1}{\rho}\Big),\quad |\rho|\to\iy,\;
\rho\in \overline{S},\;x\in [\al,T], \;\nu=\overline{0,n-1},    \eqno(1.11)
$$
where
$$
\be_1(x)=-\di\frac{1}{n}\di\int_\al^{x}p_{n-2}(t)\,dt.         \eqno(1.12)
$$
Thus, for $y_k^{(\nu)}(x,\rho),\;k=\overline{1,n},\;
\nu=\overline{0,n-1},$ we obtain the asymptotics
$$
y_k^{(\nu)}(x,\rho)=(\rho R_k)^{\nu}\exp(\rho R_k x)
\Big(1+\di\frac{\be_1(x)}{\rho R_k}+o\Big(\di\frac{1}{\rho}\Big)\Big),
\quad |\rho|\to\iy,\; \rho\in \overline{S},                    \eqno(1.13)
$$
uniformly in $x \in [\al,T].$

\medskip
Furthermore, since $z_{\nu k}(x,\rho)$ are solutions of (1.5), the functions
$y_k^{(\nu)}(x,\rho)$ satisfy (1.4). Consequently, the functions
$y_k(x,\rho)$ are solutions of the differential equation (1.1).
Thus, we arrive at the following theorem.

\bigskip
{\bf Theorem 1.1.} {\it Fix $\al\in [0,T]$ and a sector $S$ with the
property (1.2). For $\rho\in \overline{S},\; |\rho|\ge \rho_\al,\;
x \in [0,T],$ there exists a FSS $\{y_k(x,\rho)\}_{k=\overline{1,n}}$
of the equation (1.1) such that:\\
1) the functions $y_k^{(\nu)}(x,\rho),\; \nu=\overline{0,n-1}$ are continuous
for $x\in [0,T],\; \rho\in \overline{S},\; |\rho|\ge \rho_\al;$ \\
2) for each $x \in [0,T]$, the functions $y_k^{(\nu)}(x,\rho),\;
\nu=\overline{0,n-1}$ are analytic with respect to $\rho \in S,\;
|\rho|\ge \rho_\al$;\\
3) uniformly for $x\in [\al,T],$ the asymptotic formula (1.13) holds,
where $\be_1(x)$ is defined by (1.12); \\
4) as $|\rho|\to\iy,\; \rho\in \overline{S},$}
$$
\det[y_k^{(\nu-1)}(x,\rho)]_{\nu,k=\overline{1,n}}=
\rho^{\frac{n(n-1)}{2}}\det[R_k^{\nu-1}]_{\nu,k=\overline{1,n}}\Big(1+
o\Big(\di\frac{1}{\rho}\Big)\Big).
$$

\bigskip
{\bf 1.2.} In this section we discuss the possibility to obtain more precise
asymptotic formulae than (1.13). For this purpose we need some smoothness
for the coefficients of equation (1.1). Denote by $W_N[a,b]$ the set of
functions $f(x),\; x \in[a,b]$ such that the functions $f^{(\nu)}(x),\;
\nu=\overline{0,N-1}$ are absolutely continuous, and $f^{(\nu)}(x)\in L(a,b),
\;\nu=\overline{0,N}$. For $N\le 0$ we put $W_N[a,b]=L(a,b).$

Suppose that $p_{n-2}(x) \in W_1[\al,T].$
Substituting (1.11) into the right-hand side of (1.5), we get
$$
z_{\nu k}(x,\rho)=1-\di\frac{1}{n(\rho R_k)}\di\int_\al^{x}
p_{n-2}(t)\Big(1+\di\frac{\be_1(t)}{\rho R_k}\Big)\,dt-
\di\frac{1}{n(\rho R_k)^{2}}\di\int_\al^{x}p_{n-3}(t)\,dt
$$
$$
-\di\frac{1}{n\rho}\di\int_\al^{x}\Big(\di\sum_{j=1}^{k-1} R_j^{\nu +1}
R_k^{n-\nu-2}\exp(\rho (R_j-R_k)(x-t))\Big)p_{n-2}(t)\,dt
$$
$$
+\di\frac{1}{n\rho}\di\int_x^{T}\Big(\di\sum_{j=k+1}^{n}
R_j^{\nu+1}R_k^{n-\nu-2} \exp(\rho (R_j-R_k)(x-t))\Big)p_{n-2}(t)\,dt+
o\Big(\di\frac{1}{\rho^{2}}\Big).
$$
Integration by parts yields
$$
z_{\nu k}(x,\rho)=1-\di\frac{1}{n(\rho R_k)}\di\int_\al^{x}p_{n-2}(t)\,dt-
\di\frac{1}{n(\rho R_k)^{2}}\di\int_\al^{x}\Big(p_{n-3}(t)+
p_{n-2}(t)\be_1(t)\Big)\,dt
$$
$$
-\di\frac{1}{n(\rho R_k)^{2}}\di\sum_{j=1,\atop j\ne k}^{n}
\di\frac{R_j^{\nu +1} R_k^{-\nu}}{R_k-R_j}p_{n-2}(x)+
\di\frac{1}{n(\rho R_k)^{2}}\di\sum_{j=1}^{k-1}\di\frac{ R_j^{\nu +1}
R_k^{-\nu}}{R_k-R_j} \exp(\rho (R_j-R_k)(x-\al))p_{n-2}(\al)
$$
$$
+\di\frac{1}{n(\rho R_k)^{2}}
\di\sum_{j=k+1}^{n}\di\frac{ R_j^{\nu +1} R_k^{-\nu}}{R_k-R_j}
\exp(\rho (R_j-R_k)(x-T))p_{n-2}(T)+o\Big(\di\frac{1}{\rho^{2}}\Big).
$$

This asymptotic formula contains the terms with exponentials, and it is
not convenient for applications. In order to obtain the asymptotics
without terms having exponentials, we should modify slightly the
original equation (1.3). More precisely, instead of (1.3) we consider
the following equation:
$$
y_k(x,\rho)=\di\sum_{j=1}^{n}C_j\exp(\rho R_j x)- \di\frac{1}{n\rho^{n-1}}
\di\int_\al^{x}\Big(\di\sum_{j=1}^{k}R_j\exp(\rho R_j(x-t))\Big)
\Big(\di\sum_{m=0}^{n-2} p_m(t) y_k^{(m)}(t,\rho)\Big)\,dt
$$
$$
+\di\frac{1}{n\rho^{n-1}}
\di\int_x^{T}\Big(\di\sum_{j=k+1}^{n} R_j\exp(\rho R_j(x-t))\Big)
\Big(\di\sum_{m=0}^{n-2}p_m(t) y_k^{(m)}(t,\rho)\Big)\,dt,    \eqno(1.14)
$$
where
$$
C_j=
\left\{\begin{array}{ll}
-\di\frac{p_{n-2}(\al)}{n(\rho R_k)^{2}}\di\frac{R_j}{R_k-R_j}
\exp(\rho (R_k-R_j) \al),\quad & \mbox{ for } \; j=\overline{1,k-1}, \\[3mm]
1,\quad & \mbox{ for } \; j=k,\\[3mm]
-\di\frac{p_{n-2}(T)}{n(\rho R_k)^{2}}\di\frac{R_j}{R_k-R_j}
\exp(\rho (R_k-R_j) T),\quad & \mbox{ for } \; j=\overline{k+1,n}.
\end{array}\right.
$$
This implies
$$
y_k^{(\nu)}(x,\rho)=\di\sum_{j=1}^{n}C_j(\rho R_j)^{\nu}\exp(\rho R_j x)
$$
$$
-\di\frac{1}{n\rho^{n-1}}\di\int_\al^{x}
\Big(\di\sum_{j=1}^{k}R_j (\rho R_j)^{\nu}\exp(\rho R_j(x-t))\Big)
\Big(\di\sum_{m=0}^{n-2} p_m(t) y_k^{(m)}(t,\rho)\Big)\,dt
$$
$$
+\di\frac{1}{n\rho^{n-1}}\di\int_x^{T}\Big(\di\sum_{j=k+1}^{n}
R_j (\rho R_j)^{\nu}\exp(\rho R_j(x-t))\Big)
\Big(\di\sum_{m=0}^{n-2}p_m(t) y_k^{(m)}(t,\rho)\Big)\,dt.  \eqno(1.15)
$$
Denote
$$
z_{\nu k}(x,\rho)=(\rho R_k)^{-\nu}\exp(-\rho R_k x)y_k^{(\nu)}(x,\rho).
$$
Then (1.15) becomes
$$
z_{\nu k}(x,\rho)=1-\di\frac{p_{n-2}(\al)}{n(\rho R_k)^{2}}
\di\sum_{j=1}^{k-1}\di\frac{ R_j^{\nu +1} R_k^{-\nu}}{R_k-R_j}
\exp(\rho (R_j-R_k)(x-\al))
$$
$$
-\di\frac{p_{n-2}(T)}{n(\rho R_k)^{2}}\di\sum_{j=k+1}^{n}
\di\frac{ R_j^{\nu +1} R_k^{-\nu}}{R_k-R_j}\exp(\rho (R_j-R_k)(x-T))
$$
$$
-\di\frac{1}{n\rho^{n-1}}\di\int_\al^{x}
\Big(\di\sum_{j=1}^{k} R_j^{\nu+1} R_k^{-\nu}\exp(\rho (R_j-R_k)(x-t))\Big)
\Big(\di\sum_{m=0}^{n-2} p_m(t)(\rho R_k)^{m}z_{m k}(t,\rho)\Big)\,dt
$$
$$
+\di\frac{1}{n\rho^{n-1}}\di\int_x^{T}\Big(\di\sum_{j=k+1}^{n}
R_j^{\nu+1} R_k^{-\nu}\exp(\rho (R_j-R_k)(x-t))\Big)
\Big(\di\sum_{m=0}^{n-2}p_m(t) (\rho R_k)^{m}z_{m k}(t,\rho)\Big)\,dt,
$$
$$
\nu=\overline{0,n-1}.                                           \eqno(1.16)
$$
Solving the system (1.16) by the method of successive approximations
and repeating the preceding arguments, we obtain that system (1.16)
has a unique solution $z_{\nu k}(x,\rho)$ such that (1.10) holds.
Substituting (1.10) into the right-hand side of (1.16), we arrive
at (1.11). Substituting now (1.11) into the right-hand side of
(1.16) we get by the same way as above that
$$
z_{\nu k}(x,\rho)=1+\di\frac{\be_1(x)}{\rho R_k}+
\di\frac{\be_{2\nu}(x)}{(\rho R_k)^{2}}+
o\Big(\di\frac{1}{\rho^{2}}\Big),\quad
|\rho|\to\iy,\; \rho\in \overline{S},\;x \in [\al,T],       \eqno(1.17)
$$
where
$$
\be_{2\nu}(x)=-\di\frac{1}{n}\di\int_\al^{x}\Big(p_{n-3}(t)+
p_{n-2}(t)\be_1(t)\Big)\,dt-\di\frac{1}{n}\di\sum_{j=1 \atop j\ne k}^{n}
\di\frac{ R_j^{\nu +1} R_k^{-\nu}}{R_k-R_j}p_{n-2}(x).
$$
Since
$$
\di\frac{1}{n}\di\sum_{j=1 \atop j\ne k}^{n}
\di\frac{ R_j^{\nu +1} R_k^{-\nu}}{R_k-R_j}=-\di\frac{n-1}{2}+\nu,
$$
the asymptotic formula (1.17) takes the form
$$
z_{\nu k}(x,\rho)=1+\di\frac{\be_1(x)}{\rho R_k}+
\di\frac{\be_{2}(x)+\nu \be'_1(x)}{(\rho R_k)^{2}}+
o\Big(\di\frac{1}{\rho^{2}}\Big),\quad
|\rho|\to\iy,\; \rho\in \overline{S},\;x \in [\al,T],          \eqno(1.18)
$$
where
$$
\be_2(x)=-\di\frac{1}{n}\di\int_\al^{x}\Big(p_{n-3}(t)+
p_{n-2}(t)\be_1(t)\Big)\,dt-\di\frac{n-1}{2}\be'_1(x).
$$
Hence
$$
y_k^{(\nu)}(x,\rho)=(\rho R_k)^{\nu}\exp(\rho R_k x)
\Big(1+\di\frac{\be_1(x)}{\rho R_k}+
\di\frac{\be_{2}(x)+\nu \be'_1(x)}{(\rho R_k)^{2}}+
o\Big(\di\frac{1}{\rho^{2}}\Big)\Big).
$$

By the same way one can obtain the following more general assertion.

\bigskip
{\bf Theorem 1.2. }{\it Fix $\al \in [0,T],\; N\ge0$ and a sector $S$ with
the property (1.2). Assume that $p_m(x)\in W_{N+m-n+2}[\al,T],\;
m=\overline{0,n-2}.$ Then for $\rho\in\overline{S},\; |\rho|\ge \rho_\al,\;
x\in [0,T],$ there exists a FSS $\{y_k(x,\rho)\}_{k=\overline{1,n}}$ of
the equation (1.1) such that:\\
1) the functions $y_k^{(\nu)}(x,\rho),\; \nu=\overline{0,n-1},$ are
continuous for $x\in [0,T],\; \rho\in\overline{S},\; |\rho|\ge \rho_\al;$ \\
2) for each $x\in [0,T],$ the functions $y_k^{(\nu)}(x,\rho),\;\nu=
\overline{0,n-1},$ are analytic with respect to $\rho\in S,\;
|\rho|\ge \rho_\al;$ \\
3) as $|\rho|\to\iy,\; \rho\in\overline{S},$ uniformly for $x\in [\al,T],$
$$
y_k(x,\rho)=\exp(\rho R_k x)
\Big(1+\di\sum_{s=1}^{N+1}\di\frac{\be_s(x)}{(\rho R_k)^{s}}+
o\Big(\di\frac{1}{\rho^{N+1}}\Big)\Big),                        \eqno(1.19)
$$
$$
y_k^{(\nu)}(x,\rho)=(\rho R_k)^{\nu}\exp(\rho R_k x)
\Big(1+\di\sum_{s=1}^{N+1}\di\frac{\be_{s\nu}(x)}{(\rho R_k)^{s}}+
o\Big(\di\frac{1}{\rho^{N+1}}\Big)\Big),
$$
where
$$
\be_{s\nu}(x)=\di\sum_{r=0}^{\nu}C_\nu^{r}\be_{s-r}^{(r)}(x),\quad
C_\nu^{r}:=\di\frac{\nu!}{r!(\nu-r)!},
$$
$$
\be_0(x)=1,\quad \be_s(x)=0 \quad \mbox{ for }\; s<0,
$$
$$
\be^{\prime}_s(x)=-\di\frac{1}{n}\Big(\di\sum_{r=2}^{n}
C_n^{r}\be_{s+1-r}^{(r)}(x)+\di\sum_{m=0}^{n-2}p_m(x)
\be_{s-n+m+1,m}(x)\Big)\quad \mbox{\rm for}\; s\ge 1.           \eqno(1.20)
$$
In particular, (1.20) yields}
$$
\be^{\prime}_1(x)=-\di\frac{1}{n}p_{n-2}(x),
$$
$$
\be^{\prime}_2(x)=-\di\frac{1}{n}
\Big(C_n^{2}\be^{\prime\prime}_1(x)+p_{n-2}(x)\be_1(x)+p_{n-3}(x)\Big),
$$
$$
\be^{\prime}_3(x)=-\di\frac{1}{n}
\Big(C_n^{2}\be^{\prime\prime}_2(x)+C_n^{3}\be^{\prime\prime\prime}_1(x)+
p_{n-2}(x)\Big(\be_2(x)+C_{n-2}^{1}\be^{\prime}_1(x)\Big)+
p_{n-3}(x)\be_1(x)+p_{n-4}(x)\Big),
$$
$$
\ldots\ldots
$$

\bigskip
The reccurent formula (1.20) for the coefficients $\be_s(x), s\ge 1,$
can be obtained by substitution (1.19) into (1.1). We note that
$\be_s(x)\in W_{N-s+2}[\al,T].$
\bigskip

\begin{center}
{\bf II. First-order systems of differential equations}\\[0.3cm]
\end{center}

{\bf 2.1.}$\;$ Consider the following system:
$$
LY(x):=
\di\frac{1}{\rho}\,Y'(x)-A(x,\rho)Y(x),\quad 0\le x\le T\le \iy, \eqno(2.1)
$$
where $Y=[y_\nu]_{\nu=\overline{1,n}}$ is a column-vector, and the matrix
$A(x,\rho)$ has the form
$$
A(x,\rho)=A_{(0)}+\di\sum_{\mu=1}^{\iy}\di\frac{A_{(\mu)}(x)}{\rho^{\mu}},
\quad A_{(\mu)}=[a_{(\mu)\nu j}]_{\nu,j=\overline{1,n}},          \eqno(2.2)
$$
We assume that\\
($i_1$) $A_0$ is a constant matrix, and its eigenvalues $R_1,R_2,\ldots,R_n$
are such that $R_k\ne 0,\; R_j\ne R_k\;(j\ne k);$\\
($i_2$) for $\mu\ge 1,\; \nu,j=\overline{1,n},\quad a_{(\mu)\nu j}(x)\in
L(0,T),$ and $\|a_{(\mu)\nu j}(x)\|_{L(0,T)}\le Ca_*^{\mu},\; a_*\ge 0.$

\bigskip
We shall say that $L\in\Lambda_0,$ if ($i_1$)-($i_2$) hold. We also consider
the classes  $\Lambda_N \subset \Lambda_0, \; (N\ge 1)$ with additional
smoothness properties of $A_{(\mu)}(x).$ More precisely, we shall say that
$L\in\Lambda_N$ if ($i_1$)-($i_2$) hold, and $a_{(\mu)\nu j}(x)\in
W_{N-\mu+1}[0,T]$ for $\mu=\overline{1,N},\;\nu,j=\overline{1,n}.$
Note that $N=1$ is the most popular case in applications.

\bigskip
{\bf 2.2. } In Sections 2.2-2.3 we provide some formal calculations in order
to show ideas. For the explicit results see Section 2.4.

Acting formally one can seek solutions $Y_k(x,\rho),\; k=\overline{1,n},$
of system (2.1) in the form
$$
Y_k(x,\rho)=\exp(\rho R_k x)
\di\sum_{\mu=0}^{\iy}\di\frac{g_{(\mu) k}(x)}{\rho^{\mu}},      \eqno(2.3)
$$
where $g_{(\mu)k}(x)=[g_{(\mu)\nu k}(x)]_{\nu=\overline{1,n}}$ are
column-vectors. Substituting (2.3) into (2.1) we get formally
$$
LY_k(x,\rho)=\exp(\rho R_k x)
\di\sum_{\mu=0}^{\iy}\di\frac{1}{\rho^{\mu}}\Big\{
\Big(R_k g_{(\mu) k}(x)+g_{(\mu-1) k}'(x)\Big)
$$
$$
-\Big(A_{(0)}g_{(\mu) k}(x)+A_{(1)}(x)g_{(\mu-1) k}(x)+\ldots
+A_{(\mu)}(x)g_{(0) k}(x)\Big)\Big\}=0
$$
(here we put $g_{(\mu)k}(x)=0$ for $\mu<0$). This yields the following
reccurent formulae for constructing the coefficients $g_{(\mu)k}(x)$:
$$
\left. \begin{array}{c}
A_{(0)}g_{(0) k}(x)-R_k g_{(0) k}(x)=0,\\[3mm]
A_{(0)}g_{(\mu) k}(x)-R_k g_{(\mu) k}(x)=\\[3mm]
g_{(\mu-1) k}'(x)-\Big(A_{(1)}(x)g_{(\mu-1) k}(x)+\ldots+
A_{(\mu)}(x)g_{(0) k}(x)\Big),\;\mu\ge 1.
\end{array}\right\}                                              \eqno(2.4)
$$
For example, if $A_{(0)}=diag[R_k]_{k=\overline{1,n}},$ then one can take
$$
g_{(0)\nu k}(x)=Q_k(x)\delta_{\nu k},\quad
Q_k(x)=\exp\Big(\di\int_0^{x} a_{(1)kk}(\xi)\,d\xi\Big),         \eqno(2.5)
$$
where $\delta_{\nu k}$ is the Kronecker symbol. Then
$$
g_{(1)\nu k}(x)=\di\frac{a_{(1)\nu k}(x)Q_k(x)}{R_k-R_\nu},\quad \nu\ne k,
$$
$$
g_{(1)kk}'(x)=a_{(2)kk}(x)Q_k(x)
+\di\sum_{\nu=1}^{n}a_{(1)k\nu}(x)g_{(1)\nu k}(x),
$$
$$
\ldots\ldots
$$

\bigskip
{\bf 2.3. } Denote
$$
L^{*}Z(x):=\di\frac{1}{\rho}\,Z'(x)+Z(x)A(x,\rho),
$$
where $Z=[z_1,\ldots,z_n]$ is a row-vector. Then
$$
ZLY+L^{*}ZY=\di\frac{1}{\rho}\,\di\frac{d}{dx}(ZY),               \eqno(2.6)
$$
and consequently,
$$
\di\int_a^{b} ZLY\,dx=\di\frac{1}{\rho}\,(ZY)\Big|_{a}^{b}
-\di\int_a^{b} L^{*}ZY\,dx.                                       \eqno(2.7)
$$

Let $U(x)=[u_{jk}(x)]_{j,k=\overline{1,n}}$ be an absolutely continuous
non-degenerate matrix, and let $V(x)=[v_{jk}(x)]_{j,k=\overline{1,n}}$ be
the inverse matrix: $V(x)=(U(x))^{-1}.$ Denote
$$
U_k(x)=\left[ \begin{array}{l}
u_{1k}(x) \\ \ldots \\ u_{nk}(x)
\end{array}\right], \qquad
V_j(x)=[v_{j1}(x),\ldots, v_{jn}(x)].
$$

Let us show that system (2.1) is equivalent to the integral equation
$$
Y(x)=\di\sum_{j=1}^{n}U_j(x)\Big(I_j+\rho\di\int_{\ga_j}^{x}
L^{*}V_j(t)Y(t)\,dt\Big),\quad \ga_j\in [0,T],                \eqno(2.8)
$$
where $I_j$ are arbitrary constants.

Indeed, assume that $Y(x)$ satisfies (2.1), i.e. $LY(x)=0$.
Then, by virtue of (2.7),
$$
(V_j Y)\Big|_{\ga_j}^{x}=
\rho\di\int_{\ga_j}^{x}L^{*}V_j(t)Y(t)\,dt,\quad j=\overline{1,n}, \eqno(2.9)
$$
where $\ga_j\in [0,T]$ are arbitrary fixed numbers. This implies
$$
V_j(x)Y(x)=I_j+\rho\di\int_{\ga_j}^{x}L^{*}V_j(t)Y(t)\,dt,
\quad j=\overline{1,n},                                          \eqno(2.10)
$$
where $I_j=V_j(\ga_j)Y(\ga_j).$ Since $V=U^{-1},$ we arrive at (2.8).

Inversly, if Y(x) satisfies (2.8) with certain constants $I_j$ and
$\ga_j\in [0,T],$ then (2.10) is valid. In particular, this yields
$I_j=V_j(\ga_j)Y(\ga_j),$ and consequently, one gets (2.9). Taking
now (2.7) into account, we obtain
$$
\di\int_{\ga_j}^{x} V_j(t)LY(t)\,dt=0,\quad j=\overline{1,n},
$$
i.e. $V_j(x)LY(x)=0,\; j=\overline{1,n}.$ Since $\det V(x)\ne 0,$
we get $LY(x)=0,$ i.e. $Y(x)$ is a solution of the equation (2.1).

\medskip
Using the integral equation (2.8) with concrete $U(x), I_j$ and $\ga_j,$
one can obtain various solutions of system (2.1) having desirable properties.

\bigskip
{\bf 2.4. } In this section we construct the Birkhoff-type solutions for
system (2.1). The main result is formulated in Theorem 2.1.

Fix a sector $S$ in the $\rho$-plane such that
$$
Re(\rho R_1)<\ldots<Re(\rho R_n),\quad \rho\in S.                \eqno(2.11)
$$

\bigskip
{\bf Theorem 2.1. }{\it Fix $N\ge 0,$ and assume that $L\in\Lambda_N.$
Then, there exists $\rho_{*}>0$ such that for $\rho\in\overline{S},\;
|\rho|\ge \rho_{*},\; x\in [0,T],$ system (2.1) has a FSS
$\{Y_k(x,\rho)\}_{k=\overline{1,n}}$ with the following properties:\\
1) $Y_k(x,\rho)$ are continuous for $x\in [0,T],\; \rho\in\overline{S},\;
|\rho|\ge \rho_{*};$\\
2) for each $x\in [0,T],$ $Y_k(x,\rho)$ are analytic with respect to
$\rho\in S,\; |\rho|\ge \rho_{*};$\\
3) for $\rho\in\overline{S},\; |\rho|\ge \rho_{*},\;x\in [0,T],$
$$
Y_k(x,\rho)=\exp(\rho R_k x)
\Big(\di\sum_{\mu=0}^{N}\di\frac{g_{(\mu) k}(x)}{\rho^{\mu}}+
\di\frac{\ee_k(x,\rho)}{\rho^{N}}\Big),                        \eqno(1.12)
$$
where the column-vectors $g_{(\mu) k}(x)=[g_{(\mu) \nu k}(x)]_{\nu=
\overline{1,n}}$ are defined by (2.4), and for the vector $\ee_k(x,\rho)=
[\ee_{\nu k}(x,\rho)]_{\nu=\overline{1,n}},$
$$
\di\lim_{|\rho|\to\iy\atop\rho\in\overline{S}}\di\max_{0\le x\le T}|
\ee_{\nu k}(x,\rho)|=0,\quad \nu,k=\overline{1,n},
$$
i.e. uniformly in $x\in [0,T],\; \ee_k(x,\rho)=o(1)$ as $|\rho|\to\iy,
\;\rho\in\overline{S}.$}

\bigskip
{\it Proof.} Without loss of generality we consider here the case
when $A_{(0)}$ is a diagonal matrix:
$$
A_{(0)}=diag[R_k]_{k=\overline{1,n}}.
$$
The general case is studied in the end of the proof by reduction
to the diagonal case.

Define the matrix $U(x,\rho)=[u_{jk}(x,\rho)]_{j,k=\overline{1,n}}$
with the columns $U_k(x,\rho)=[u_{jk}(x,\rho)]_{j=\overline{1,n}}$
by the formula
$$
U_k(x,\rho)=\exp(\rho R_k x)\di\sum_{\mu=0}^{N}
\di\frac{g_{(\mu) k}(x)}{\rho^{\mu}},\quad k=\overline{1,n},   \eqno(2.13)
$$
where $g_{(\mu) k}(x)$ satisfy (2.4) and (2.5) holds. Let $V(x,\rho)=
(U(x,\rho))^{-1}$ be inverse matrix with the rows $V_j(x,\rho)=
[v_{j1}(x,\rho),\ldots, v_{jn}(x,\rho)].$ For each fixed $k=\overline{1,n}$
and $\rho\in\overline{S},$ we consider the integral equation
$$
Y_k(x,\rho)=U_k(x,\rho)+\rho\di\int_0^{x}
\Big(\di\sum_{j=1}^{k}U_j(x,\rho)L^{*}V_j(t,\rho)\Big)Y_k(t,\rho)\,dt
$$
$$
-\rho\di\int_x^{T}\Big(\di\sum_{j=k+1}^{n}
U_j(x,\rho)L^{*}V_j(t,\rho)\Big)Y_k(t,\rho)\,dt,              \eqno(2.14)
$$
with respect to the column-vector $Y_k.$ We note that (2.14) is a
particular case of (2.8) when
$$
I_j=
\left\{ \begin{array}{ll}
0,\quad & j\ne k, \\
1,\quad & j=k,
\end{array}\right. \qquad
\ga_j=
\left\{ \begin{array}{ll}
0,\quad & j\le k,\\
T,\quad & j>k.
\end{array}\right.
$$

In order to study the solvability of (2.14) we need some preliminary
calculations. Denote by $\Gamma_a$ the set of functions $\theta(x,\rho)$
of the form
$$
\theta(x,\rho)=\di\sum_{\mu=0}^{\iy}\di\frac{\theta_{(\mu)}(x)}{\rho^{\mu}},
$$
where $\theta_{(\mu)}(x)\in L(0,T),$ and $\|\theta_{(\mu)}(x)\|_{L(0,T)}
\le Ca^{\mu}.$

It follows from (2.13), (2.4) and (2.5) that
$$
LU_k(x,\rho)=\di\frac{1}{\rho^{N+1}}\exp(\rho R_k x)\Big(
H_{(0) k}(x)+\di\frac{H_{(1) k}(x,\rho)}{\rho}\Big),           \eqno(2.15)
$$
where $H_{(\mu) k}=[H_{(\mu)\nu k}]_{\nu=\overline{1,n}}$ are
column-vectors such that
$$
H_{(0) k}(x)=g_{(N) k}'(x)-
\Big(A_{(1)}(x)g_{(N) k}(x)+\ldots+A_{(N+1)}(x)g_{(0) k}(x)\Big),\eqno(2.16)
$$
$$
H_{(1)\nu k}(x,\rho)\in\Gamma_{a_*}.                             \eqno(2.17)
$$
In particular, (2.16) yields
$$
H_{(0) kk}(x)=0,\quad H_{(0)\nu k}(x)\in L[0,T],\; \nu\ne k.     \eqno(2.18)
$$

The next step is to calculate the inverse matrix $V(x,\rho)=
(U(x,\rho))^{-1}.$ Since
$$
g_{(0) k}(x)=[\delta_{\nu k} Q_k(x)]_{\nu=\overline{1,n}},
$$
it follows from (2.13) that
$$
\Big(\det U(x,\rho)\Big)^{-1}= \exp\Big(-\rho x \di\sum_{k=1}^{n} R_k\Big)
\Big(\di\prod_{k=1}^{n}Q_k(x)\Big)^{-1}\Big(1+\di\frac{\theta(x,\rho)}{\rho}
\Big),
$$
and $\theta(x,\rho)\in\Gamma_{a_1}$ for a certain $a_1>0.$ Therefore,
$$
V_j(x,\rho)=\exp(-\rho R_j x)\Big(g_{(0) j}^{*}(x)+
\di\frac{g_{(1)j}^{*}(x,\rho)}{\rho}\Big),\quad j=\overline{1,n}, \eqno(2.19)
$$
where $g_{(\mu)j}^{*}=[g_{(\mu)j1}^{*},\ldots,g_{(\mu)jn}^{*}],\; \mu=0,1,$
are row-vectors such that
$$
g_{(0)j\nu}^{*}(x)=\Big(Q_j(x)\Big)^{-1}\delta_{j \nu},
\quad g_{(1)j\nu}^{*}(x,\rho)\in\Gamma_{a_1}.                   \eqno(2.20)
$$

By virtue of (2.15), (2.17)-(2.20),
$$
VLU=\Big[\di\frac{1}{\rho^{N+1}}\Big(h_{(0)jk}(x)+
\di\frac{h_{(1)jk}(x,\rho)}{\rho}\Big)
\exp(\rho(R_k-R_j)x)\Big]_{j,k=\overline{1,n}},                \eqno(2.21)
$$
where
$$
h_{(0)jk}(x)=\Big(Q_j(x)\Big)^{-1}H_{(0)jk}(x),\; h_{(1)jk}(x,\rho)\in
\Gamma_{a_2},\quad j,k=\overline{1,n},\; a_2=\max(a_*, a_1).    \eqno(2.22)
$$
In particular, $h_{(0)kk}(x)=0.$

Furthermore, since $VU=E$ is the indentity matrix, it follows from (2.6) that
$$
VLU+L^{*}VU=0,
$$
i.e.
$$
L^{*}V=-VLUV.                                                    \eqno(2.23)
$$
Using (2.23) and (2.19)-(2.22) we calculate
$$
L^{*}V_j(x,\rho)=\di\frac{1}{\rho^{N+1}}\exp(-\rho R_j x)
\Big(\omega_{(0)j}(x)+\di\frac{\omega_{(1)j}(x,\rho)}{\rho}\Big),
\quad j=\overline{1,n},                                         \eqno(2.24)
$$
where $\omega_{(\mu)j}=[\omega_{(\mu)j1},\ldots,\omega_{(\mu)jn}],
\;\mu=0,1,$ are row-vectors such that
$$
\omega_{(0)jk}(x)=-\Big(Q_j(x) Q_k(x)\Big)^{-1}H_{(0)jk}(x),\;
\omega_{(1)jk}(x,\rho)\in\Gamma_{a_2},\quad j,k=\overline{1,n}. \eqno(2.25)
$$
In particular,
$$
\omega_{(0)kk}(x)=0,\quad k=\overline{1,n}.                     \eqno(2.26)
$$

Denote
$$
W_k^{0}(x,\rho)=
\di\sum_{\mu=0}^{N}\di\frac{g_{(\mu) k}(x)}{\rho^{\mu}}.        \eqno(2.27)
$$
By the replacement
$$
Y_k(x,\rho)=\exp(\rho R_k x)W_k(x,\rho),                        \eqno(2.28)
$$
we transform the integral equation (2.14) to the form
$$
W_k(x,\rho)=W_k^{0}(x,\rho)+
$$
$$
\di\frac{1}{\rho^{N}}\di\int_0^{x}
\Big(\di\sum_{j=1}^{k}\exp\Big(\rho (R_j-R_k) (x-t)\Big)
W_j^{0}(x,\rho)\Big(\omega_{(0)j}(t)+
\di\frac{\omega_{(1)j}(t,\rho)}{\rho}\Big)\Big)W_k(t,\rho)\,dt
$$
$$
-\di\frac{1}{\rho^{N}}\di\int_x^{T}\Big(\di\sum_{j=k+1}^{n}
\exp\Big(\rho (R_j-R_k) (x-t)\Big)W_j^{0}(x,\rho)\Big(\omega_{(0)j}(t)
+\di\frac{\omega_{(1) j}(t,\rho)}{\rho}\Big)\Big)W_k(t,\rho)\,dt, \eqno(2.29)
$$
or
$$
W_k(x,\rho)=W_k^{0}(x,\rho)
+\di\frac{1}{\rho^{N}}\di\int_0^{T} B_k(x,t,\rho)W_k(t,\rho)\,dt, \eqno(2.30)
$$
where the matrix $B_k=[B_{k\nu s}]_{\nu,s=\overline{1,n}}$ has the form
$$
B_k(x,t,\rho)=
\left\{ \begin{array}{ll}
\di\sum_{j=1}^{k}\exp\Big(\rho (R_j-R_k) (x-t)\Big)W_j^{0}(x,\rho)
\Big(\omega_{(0) j}(t)+
\di\frac{\omega_{(1) j}(t,\rho)}{\rho}\Big),\, & x\ge t,\\
-\di\sum_{j=k+1}^{n}\exp\Big(\rho (R_j-R_k) (x-t)\Big)W_j^{0}(x,\rho)
\Big(\omega_{(0) j}(t)+\di\frac{\omega_{(1) j}(t,\rho)}{\rho}\Big),\, & x<t.
\end{array}\right.                                                \eqno(2.31)
$$
Solving (2.30) by the method of successive approximations and using
(2.24)-(2.27), we obtain that there exists $\rho_{*}>0$ such that for
$\rho\in\overline{S},\; |\rho|\ge\rho_{*},\; x\in [0,T],\; k=
\overline{1,n},$ the integral equation (2.30) has a unique solution
$W_k(x,\rho)$ having the following asymptotics
$$
W_k(x,\rho)=W_k^{0}(x,\rho)+o\Big(\di\frac{1}{\rho^{N}}\Big),
\quad |\rho|\to\iy,\;\rho\in\overline{S},                        \eqno(2.32)
$$
uniformly for $x\in [0,T].$

Indeed, let $N\ge 1.$ It follows from (2.31), (2.11), (2.25) and (2.27) that
$$
\di\max_{0\le x\le T}\di\int_0^{T} |B_{k\nu s}(x,t,\rho)|\,dt\le C,\quad
|\rho|\ge a_3,\; \rho\in\overline{S},\; k,\nu,s=\overline{1,n},\; a_3>a_2.
$$
Then the method of successive approximations gives
$$
W_k(x,\rho)=W_k^{0}(x,\rho)+O\Big(\di\frac{1}{\rho^{N}}\Big),\quad
|\rho|\to\iy,\;\rho\in\overline{S},\; k=\overline{1,n},          \eqno(2.33)
$$
uniformly in $x\in [0,T].$ Substituting (2.33) into the right-hand side
of (2.29) we get
$$
W_k(x,\rho)=W_k^{0}(x,\rho)+\di\frac{1}{\rho^{N}}\di\int_0^{x}
\Big(\di\sum_{j=1}^{k}\exp\Big(\rho (R_j-R_k) (x-t)\Big)g_{(0)j}(x)
\omega_{(0) j}(t)g_{(0)k}(t)\,dt
$$
$$
-\di\frac{1}{\rho^{N}}\di\int_x^{T}\Big(\di\sum_{j=k+1}^{n}
\exp\Big(\rho (R_j-R_k) (x-t)\Big)g_{(0)j}(x)\omega_{(0) j}(t)g_{(0)k}(t)\,dt
+O\Big(\di\frac{1}{\rho^{N+1}}\Big).
$$
The terms with $j\ne k$ give us $o\Big(\di\frac{1}{\rho^{N}}\Big)$ as
$\rho\to\iy,\; \rho\in\overline{S}$ uniformly in $x\in [0,T].$ The term
with $j=k$ is equal to zero, since according to (2.5) and (2.26),
$$
\omega_{(0)k}(t)g_{(0)k}(t)=0.
$$
Thus, for $N\ge 1,$ (2.32) is proved. The case $N=0$ requires different
calculations (see [7], [9] for details).

\medskip
Using (2.28), (2.32) and (2.27), we arrive at (2.12), and consequently
Theorem 2.1 is proved for the case when $A_{(0)}$ is a diagonal matrix.

\bigskip
Now we study the general case when $A_{(0)}$ is an arbitrary matrix with
the property ($i_1$). Let
$$
\Omega_k=[\Omega_{\nu k}]_{\nu=\overline{1,n}},\quad k=\overline{1,n},
$$
be eigenvectors of $A_{(0)}$ for the eigenvalues $R_k,\; k=\overline{1,n}.$
Consider the matrix
$$
\Omega=[\Omega_1,\ldots,\Omega_n]=[\Omega_{\nu k}]_{\nu, k=\overline{1,n}}.
$$
By the replacement
$$
Y(x)=\Omega\tilde Y(x),
$$
we transform system (2.1) to the form
$$
\di\frac{1}{\rho}\,\tilde Y'(x)=\tilde A(x,\rho)\tilde Y(x),   \eqno(2.34)
$$
where
$$
\tilde A(x,\rho)=\Omega^{-1} A(x,\rho)\Omega.
$$
Clearly,
$$
\tilde A(x,\rho)=\tilde A_{(0)}+\di\sum_{\mu=0}^{\iy}
\di\frac{\tilde A_{(\mu)}(x)}{\rho^{\mu}},\quad
\tilde A_{(\mu)}=[\tilde a_{(\mu)\nu j}]_{\nu,j=\overline{1,n}},
$$
and
$$
\tilde A_{(\mu)}=\Omega^{-1} A_{(\mu)}\Omega,\quad
\tilde A_{(0)}=diag[R_k]_{k=\overline{1,n}}.
$$
For system (2.34), Theorem 2.1 has been already proved. Thus, Theorem 2.1
holds also for an arbitrary $L\in\Lambda_N.$ Moreover,
$$
g_{(0) k}(x)=\tilde Q_k(x)\Omega_k,\quad
\tilde Q_k(x)=\exp\Big(\di\int_0^{x}\tilde a_{(1) kk}(\xi)\,d\xi\Big).
$$

\bigskip
{\bf 2.5. } In order to prove Theorem 2.1 one can also use the following
arguments.

Fix a sector $S$ with the property (2.11) and consider system (2.1) for
$\rho\in\overline{S}$:
$$
\di\frac{1}{\rho}\,Y'=A(x,\rho)Y                               \eqno(2.35)
$$
with respect to the matrix $Y=[y_{jk}]_{j,k=\overline{1,n}}.$ Let for
definiteness, $A_{(0)}=diag[R_k]_{k=\overline{1,n}}.$

By the replacement
$$
Y=U\xi,\quad \xi=[\xi_{jk}]_{j,k=\overline{1,n}},
$$
where the matrix $U$ is defined by (2.13), we reduce (2.35) to the system
$$
\di\frac{1}{\rho}\,\xi'=-VLU\xi                                 \eqno(2.36)
$$
with respect to $\xi.$ For the matrix $VLU$ we have the representation
(2.21). Then (2.36) becomes
$$
\xi_{jk}^{\prime}=-\di\frac{1}{\rho^{N}}\di\sum_{\nu=1}^{n}
h_{j \nu}(x,\rho)\exp\Big(\rho (R_\nu-R_j) x\Big)\xi_{\nu k},
\quad j,k=\overline{1,n},                                       \eqno(2.37)
$$
where
$$
h_{j \nu}(x,\rho)=h_{(0)j \nu}(x)+\di\frac{h_{(1)j \nu}(x,\rho)}{\rho}.
$$
Consider the integral equations
$$
\xi_{j k}(x,\rho)=\delta_{j k}-\di\frac{1}{\rho^{N}}
\di\int_{\ga_{j k}}^{x}\di\sum_{\nu=1}^{n} h_{j \nu}(t,\rho)
\exp\Big(\rho (R_\nu-R_j)t\Big)\xi_{\nu k}(t,\rho)\,dt,\;
\ga_{jk}=\left\{ \begin{array}{ll}
0,\; & j\le k,\\ T,\; & j>k.
\end{array}\right.                                               \eqno(2.38)
$$
Clearly, if $\{\xi_{j k}\}$ is a solution of (2.38), then $\{\xi_{j k}\}$
satisfy (2.37).

By the replacement
$$
\xi_{j k}(x,\rho)=\exp\Big(\rho (R_k-R_j) x\Big)\eta_{j k}(x,\rho),
$$
we reduce (2.38) to the system
$$
\eta_{j k}(x,\rho)=\delta_{j k}-\di\frac{1}{\rho^{N}}\di\int_{\ga_{j k}}^{x}
\exp\Big(\rho (R_j-R_k) (x-t)\Big)\di\sum_{\nu=1}^{n}h_{j \nu}(t,\rho)
\eta_{\nu k}(t,\rho)\,dt,\; j,k=\overline{1,n}.                   \eqno(2.39)
$$
Solving (2.39) by the method of successive approximations we obtain that
there exists $\rho_{*}>0$ such that for $\rho\in\overline{S},\; |\rho|\ge
\rho_{*},\; x\in [0,T],$ system (2.39) has a unique solution
$\eta_{j k}(x,\rho)$ having the asymtotics
$$
\eta_{jk}(x,\rho)=\delta_{jk}+o\Big(\di\frac{1}{\rho^{N}}\Big),
\rho\to\iy,\;\rho\in\overline{S},                                \eqno(2.40)
$$
uniformly in $x\in [0,T].$
Sinse
$$
y_{jk}(x,\rho)=\di\sum_{\nu=1}^{n}u_{j\nu}(x,\rho)\eta_{\nu k}(x,\rho)
\exp\Big(\rho (R_k-R_\nu)x\Big),
$$
we infer for the columns $Y_k=[y_{jk}]_{j,k=\overline{1,n}}$:
$$
Y_k(x,\rho)=\exp(\rho R_k x)\di\sum_{\nu=1}^{n}W_{\nu}^{0}(x,\rho)
\eta_{\nu k}(x,\rho).
$$
Taking (2.40) into account we arrive at the assertions of Theorem 2.1.

\bigskip
{\bf 2.8. } In this section we consider the n-th order differential
equation of the form
$$
\ell y:=y^{(n)}+\di\sum_{k=0}^{n-1}P_k(x,\rho)y^{(k)}=0,
\quad 0\le x\le T\le\iy,                                        \eqno(2.41)
$$
where
$$
P_k(x,\rho)=\rho^{n-k}p_{kk}+\rho^{n-k+1}p_{k,k+1}(x)+\ldots +p_{kn}(x).
$$
We assume that\\
($i_1$) $p_{kk},\; k=\overline{0,n-1}$ are constants, $p_{00}\ne 0,$ and
the roots $R_1,\ldots,R_n$ of the characteristic polynomial
$$
F(R)=\di\sum_{k=0}^{n}p_{k k}R^{k},\quad p_{nn}:=1,
$$
are simple;\\
($i_2$) $p_{k,k+j}(x)\in L(0,T),\; k=\overline{0,n-1},\; j\ge 1.$

We shall say that $\ell\in\Lambda_0^{\prime}$ if ($i_1$)-($i_2$) hold. If
additionally, for a certain $N\ge 1,\; p_{k,k+j}(x)\in W_{N-j+1}[0,T],\;
j=\overline{1,N},$ we shall say that $\ell\in\Lambda_N^{\prime}.$

Fix $N\ge 0$ and a sector $S$ in the $\rho$-plane with the property (2.11).

\bigskip
{\bf Theorem 2.2. }{\it Assume that $\ell\in\Lambda_N^{\prime}.$ Then there
exists $\rho_{*}>0$ such that for $\rho\in\overline{S},\; |\rho|\ge
\rho_{*},\; x\in [0,T],$ equation (2.41) has a FSS $\{y_k(x,\rho)\}_{k=
\overline{1,n}}$ with the following properties:\\
1) the functions $y_k^{(\nu-1)}(x,\rho),\; k,\nu=\overline{1,n}$ are
continuous for $x\in [0,T],\; \rho\in\overline{S},\; |\rho|\ge \rho_{*};$\\
2) for each $x\in [0,T],$ the functions $y_k^{(\nu-1)}(x,\rho),\;
k,\nu=\overline{1,n},$ are analytic with respect to $\rho\in S,\;
|\rho|\ge \rho_{*};$\\
3) as $|\rho|\to\iy,\;\rho\in\overline{S},$ uniformly in $x\in [0,T],$
$$
y_k^{(\nu-1)}(x,\rho)=(\rho R_k)^{\nu-1}\exp(\rho R_k x)
\Big(\di\sum_{\mu=0}^{N}\di\frac{G_{(\mu) \nu k}(x)}{\rho^{\mu}}
+o\Big(\di\frac{1}{\rho^{N}}\Big)\Big),\quad k,\nu=\overline{1,n},\eqno(2.42)
$$
where
$$
G_{(\mu) \nu k}(x)=\di\sum_{j=0}^{\nu-1}C_{\nu-1}^{j}R_k^{\nu-1-j}
G_{(\mu-j)0k}^{(j)},
$$
and $G_{(\mu)0k}(x)$ can be calculated by substitution (2.42) into
(2.41). In particular,}
$$
G_{(0)\nu k}=\exp\Big(\di\int_{0}^{x}\omega_k(\xi)\,d\xi\Big),
\quad \omega_k(x)=
-\di\frac{1}{F'(R_k)}\di\sum_{j=0}^{n-1}p_{j,j+1}(x)R_k^{j}.  \eqno(2.43)
$$

\bigskip
{\it Proof. } We transform (2.41) to the form
$$
\di\frac{y^{(n)}}{\rho^{n}}+\di\sum_{k=0}^{n-1}
{\cal P}_k(x,\rho)\di\frac{y^{(k)}}{\rho^{k}}=0,               \eqno(2.44)
$$
where
$$
{\cal P}_k(x,\rho)=\di\frac{1}{\rho^{n-k}}P_k(x,\rho)
=\di\sum_{\mu =0}^{n-k}\di\frac{p_{k,k+\mu}(x)}{\rho^{\mu}}.
$$
Denote
$$
y_1=y,\; y_2=\di\frac{1}{\rho}\,y',\ldots, y_n=\di\frac{1}{\rho}\,y^{(n-1)}.
$$
Then equation (2.44) is equivalent to the system
$$
\di\frac{1}{\rho}\,
\left[ \begin{array}{l}
y_1^{\prime}\\ y_2^{\prime}\\ \ldots \\ y_{n-1}^{\prime} \\ y_n^{\prime}
\end{array}\right]
=\left[ \begin{array}{ccccc}
0 & 1 & 0 & \ldots & 0 \\
0 & 0 & 1 & \ldots & 0 \\
\ldots & \ldots & \ldots & \ldots & \ldots \\
0 & 0 & 0 & \ldots & 1 \\
-{\cal P}_0(x,\rho)& -{\cal P}_1(x,\rho) & -{\cal P}_2(x,\rho) &\ldots &
-{\cal P}_{n-1}(x,\rho)
\end{array}\right]
\left[ \begin{array}{l}
y_1 \\ y_2 \\ \ldots \\ y_{n-1} \\ y_n
\end{array}\right].                                             \eqno(2.45)
$$
System (2.45) is a particular case of (2.1) with
$$
A(x,\rho)=A_{(0)}+\di\sum_{\mu=1}^{n}\di\frac{A_{(\mu)}(x)}{\rho^{\mu}},
$$
$$
A_{(0)}=
\left[ \begin{array}{ccccc}
0 & 1 & 0 & \ldots & 0 \\
0 & 0 & 1 & \ldots & 0 \\
\ldots & \ldots & \ldots & \ldots & \ldots \\
0 & 0 & 0 & \ldots & 1 \\
-p_{00} & -p_{11} & -p_{22} &\ldots &-p_{n-1,n-1}
\end{array}\right],
$$
$$
A_{(1)}(x)=
\left[ \begin{array}{cccc}
0 & 0 & \ldots & 0\\
0 & 0 & \ldots & 0\\
\ldots & \ldots & \ldots & \ldots\\
0 & 0 & \ldots & 0\\
-p_{0 1}(x) & -p_{1 2}(x) &\ldots &-p_{n-1, n}(x)
\end{array}\right],
$$
$$
A_{(2)}(x)=
\left[ \begin{array}{ccccc}
0 & 0 & \ldots & 0 & 0\\
0 & 0 & \ldots & 0 & 0\\
\ldots & \ldots & \ldots & \ldots & \ldots \\
0 & 0 & \ldots & 0 & 0\\
-p_{0 2}(x) & -p_{1 3}(x) &\ldots &-p_{n-2, n}(x) & 0
\end{array}\right],
$$
$$
A_{(n)}(x)=
\left[ \begin{array}{cccc}
0 & 0 & \ldots & 0\\
0 & 0 & \ldots & 0\\
\ldots & \ldots & \ldots & \ldots\\
0 & 0 & \ldots & 0\\
-p_{0 n}(x) & 0 &\ldots & 0
\end{array}\right].
$$
Clearly, $\Omega_k=[R_k^{\nu-1}]_{\nu=\overline{1,n}},\; k=\overline{1,n}$
are eigenvectors of $A_{(0)}$ for the eigenvalues $R_k.$ Thus, Theorem 2.2
follows from Theorem 2.1 with $\tilde a_{(1) k k}(x)=\omega_k(x),\quad
\tilde Q_k(x)=G_k(x).$
$\hfill\Box$


\begin{center}
{\bf References}
\end{center}
\begin{enumerate}
\item[{[1]}] Naimark M.A., Linear Differential Operators, 2nd ed., Nauka,
     Moscow, 1969; English transl. of 1st ed., Parts I,II, Ungar,
     New York, 1967, 1968.
\item[{[2]}] Beals R., Deift P. and Tomei C.,  Direct and Inverse Scattering
     on the Line, Math. Surveys and Monographs, v.28. Amer. Math. Soc.
     Providence: RI, 1988.
\item[{[3]}] Yurko V.A., Inverse Spectral Problems for Differential
     Operators and their Applications, Gordon and Breach, New York, 1999.
\item[{[4]}] Tamarkin J.D., On Certain General Problems of Theory of
       Ordinary Linear Differential Equations, Petrograd, 1917.
\item[{[5]}] Rasulov M.L., Contour Integral Method, Nauka, Moscow, 1964.
     English transl. Amsterdam, 1967.
\item[{[6]}] Vagabov A.I., Asymptotic behavior of solutions of differential
       equations with respect to a parameter, and applications.
       Dokl. Akad. Nauk 326 (1992), no.2, 219-223; transl. in Russian
       Acad. Sci. Dokl. Math. 46 (1993), no.2, 240-244.
\item[{[7]}] Vagabov A.I., On sharpening an asymptotic theorem of Tamarkin.
       Diff. Uravneniya 29 (1993), no.1, 41-49; transl. in Diff.
       Equations 29 (1993), no.1, 33-41.
\item[{[8]}] Rykhlov V., Asymptotics of a system of solutions of a
       differential equation of general form with a parameter. Ukrain.
       Mat. Zh. 48 (1996), no. 1, 96-108.
\item[{[9]}] Rykhlov V., Asymptotical formulas for solutions of linear
     differential systems of the first order, Results Math. 36 (1999), 342-353.
\item[{[10]}] Savchuk A.M., Shkalikov A.A. Asymptotic analysis of solutions of ordinary differential 
equations with distribution coefficients, Sb. Math. 211 (2020), no.~11, 1623-1659.

\end{enumerate}

\end{document}